\documentclass[11pt,letter]{amsart}

\usepackage{amsfonts, amsmath, amssymb, amscd, amsthm, graphicx}

\newtheorem{thm}{Theorem}[section]
\newtheorem{cor}[thm]{Corollary}
\newtheorem{lem}[thm]{Lemma}

\newtheorem{que}[thm]{Question}
\theoremstyle{definition}

\theoremstyle{remark}

\begin{document}

\title{Uniform Kazhdan groups}

\author{D. Osin}
\address{Denis Osin: Department of Mathematics, The City College of New York, New York,
NY 10031} \curraddr{} \email{denis.osin@gmail.com}
\thanks{The first author was supported by the NSF grant DMS \# 0605093}

\author{D. Sonkin}
\address{Dmitriy Sonkin: Department of Mathematics, University of Virginia, Charlottesville, VA 22904}

\subjclass[2000]{20F65, 20F67} \keywords{Kazhdan property (T),
hyperbolic group, uniform Kazhdan constant}
\email{ds5nd@virginia.edu}

\maketitle

\begin{abstract}
We construct first examples of infinite groups having property (T)
whose Kazhdan constants admit a lower bound independent of the
choice of a finite generating set.
\end{abstract}

\section{Introduction}

Let $G$ be a group generated by a finite set $X$. We say that $G$
has property (U) (with respect to $X$) if there is a constant
$C=C(G)$ such that the following condition holds. For any
generating set $Y$ of $G$, there exists an element $t\in G$ such
that $$\max\limits_{x\in X} |t^{-1}xt|_Y \le C.$$

Observe that if $G$ has property (U) with respect to $X$, it has
(U) with respect to any other finite generating set. Finite groups
are obvious examples of groups with property (U). On the other
hand, the question of the existence of infinite groups with this
property is not trivial. In this paper we prove the following.

\begin{thm}\label{U}
Every non--elementary torsion--free hyperbolic group has an
infinite quotient group with property (U).
\end{thm}

Our construction of groups with property (U) was inspired by a
well--known question about Kazhdan groups. We recall that to each
finite subset $S$ of a (discrete) group $G$ and each unitary
representation $\pi \colon G\to U(\mathcal H)$ of $G$ on a Hilbert
space $\mathcal H$, one associates the {\it Kazhdan constant}
$$\varkappa (G, S,\pi )=\inf\limits_{\xi\in \mathcal H^1}\max\limits_{s\in S } \| \pi (s)\xi -\xi \|,$$
where $\mathcal H^1$ denotes the unit sphere in $\mathcal H$.
Further we set $\varkappa(G,S)=\inf\limits_{\pi } \varkappa (G,
S,\pi )$ and $\varkappa (G)=\inf\limits_{S} \varkappa (G,S)$,
where the infimum is taken over all unitary representations $\pi $
and all finite generating sets $S$ of $G$ respectively. Recall
that a group $G$ has property (T) of Kazhdan if $\varkappa
(G,S)>0$ for some (or, equivalently, for any) generating set $S$
of $G$.

Clearly $\varkappa(G)>0$ whenever $G$ is a finite group. In
\cite{Lub}, Lubotzky asked whether this inequality holds for all
Kazhdan groups $G$. The negative answer was obtained by Gelander
and \. Zuk \cite{GZ}. They proved that any dense subgroup of a
connected locally compact topological group has zero uniform
Kazhdan constant. Shortly later Osin showed that $\varkappa(H)=0$
for any infinite hyperbolic group $H$. Thus the equality
$\varkappa(H)=0$ holds for the majority of known groups with
property (T). Moreover, no examples of infinite groups with
non--zero uniform Kazhdan constant were known until now.

It is easy to show that every group having properties (U) and (T)
has non--zero uniform Kazhdan constant. Applying Theorem \ref{U}
to a torsion--free Kazhdan hyperbolic group, we obtain the
following.

\begin{cor} \label{ut}
There exists an infinite finitely generated group $G$ such that
$\varkappa(G)>0$.
\end{cor}

The proof of Theorem \ref{U} is based on a variant of
Ol'shanskii's technique \cite{Ols} as elaborated by Semenov
\cite{Sem-diss}. It is worth to note that our group $G$ satisfies
the identity $x^n=1$ for some large odd $n$. In particular, $G$ is
not residually finite according to the positive solution of the
restricted Burnside problem \cite{Zelm1}. Our method does not
allow to avoid this identity. Thus the following question is still
open.

\begin{que}
Does there exist an infinite residually finite group $G$ such that
$\varkappa (G)$ is positive?
\end{que}

\section{The Ol'shanskii--Semenov construction}

Let us give one of many equivalent definitions of a hyperbolic
group \cite{Gro87}. A group $G$ with a finite generating set $X$
is hyperbolic (in the sense of Gromov) if its Cayley graph
$\Gamma=\Gamma(G,X)$ is a hyperbolic metric space with respect to
the word-length metric. This means that there exists a constant
$\delta$ such that every geodesic triangle in $\Gamma$ is
$\delta$-thin, i.e., each of its sides belongs to the closed
$\delta$-neighborhood of the union of the other two sides. A group
is called elementary if it contains a cyclic subgroup of finite
index.

In \cite{Ols} Ol'shanskii showed that any non-elementary torsion
free hyperbolic group has an infinite quotient satisfying the
identity $x^n=1$ provided $n$ is large and odd (later Ivanov and
Ol'shanskii \cite{IO} showed how to construct such quotients of
hyperbolic groups with torsion). The proof was based on the graded
diagrams method \cite{Ols-book}. This method was also used in
\cite{Ols_cycl}, where Ol'shanskii constructed infinite finitely
generated groups all of whose proper subgroups are finite cyclic.
In \cite{Sem-diss, Sem-obz}, Semenov incorporated the two
constructions to obtain infinite quotients of non-elementary
torsion free hyperbolic groups all of whose proper subgroups are
finite cyclic.

In this paper we use the methods of \cite{Ols-book} as elaborated
in \cite{Sem-diss, Sem-obz}. Our proofs heavily depend on
technical lemmas from \cite{Sem-diss} that are collected below
(otherwise our paper would be unreasonably long). Let us recall
the main steps of the Ol'shanskii--Semenov construction.

Let $$G=\langle a_1,~a_2,\dots ,~a_m | R=1, ~R \in \mathcal R_0
\rangle $$ be a non-elementary torsion free hyperbolic group. In
what follows, certain parameters ($C$, $d$, $h$, $n$, $n_0$)
appear. In fact, we do not need the exact values of these
parameters. For our goals it suffices to know that {\it there
exist} parameters $C$, $d$, $h$, $n$, $n_0$ such that all results
listed below are true and, in addition, $n_0 \gg n \gg 1$ and
$n_0$ is an odd number. (For exact values of these parameters and
other details we refer to \cite{Ols,Sem-diss,Sem-obz}.)

We construct a sequence of quotient group $G(i)$ of $G$ as
follows. Let $G(0)=G$. Assuming that the group $G(i-1)=\langle
a_1,~a_2,\dots ,~a_m | R=1, ~R \in \mathcal R_{i-1} \rangle$ is
already constructed, we will introduce the set $\mathcal S_i$ of
defining words of rank $i$, and set $\mathcal R_i = \mathcal
R_{i-1} \cup \mathcal S_i$ and
$$
G(i)=\langle a_1,~a_2,\dots ,~a_m | R=1, ~R \in \mathcal R_i
\rangle
$$

Elements of $G(i)$ are referred to as words over the alphabet
$\{a_1^{\pm 1}, \dots ,a_m^{\pm 1}\}$. By $|W|$ we denote the
length of a word $W$ with respect to $\{a_1^{\pm 1}, \dots
,a_m^{\pm 1}\}$. We write $A \equiv B$ to indicate
letter-for-letter equality of words $A$ and $B$.

We say that a word $W$ is simple in $G(0)$ and in all ranks $i<C$
if it is not equal to a proper power of some element of $G(0)$ and
is not conjugate to any word of smaller length. A word $W$ is
called simple in rank $i \ge C$ if it is not conjugated to a power
of a shorter word in the group $G(i-1)$ and is not conjugated to a
power of a period of rank $k \le i-1$ in the group $G(i-1)$. For
each $i=1,~2, \dots$ choose some set $\mathcal X_i$ of simple in
rank $i-1$ words of length $i$, maximal with respect to the
property that if $A,~B \in \mathcal X_i$, $A \not \equiv B$, then
$A$ is not conjugated to $B^{\pm 1}$ in the group $G(i-1)$. Words
from $\mathcal X_i$ are called periods of rank $i$.
For every period $A \in \mathcal X_i$ choose a maximal subset of
words $\mathcal G_A$ such that:
\begin{itemize}
\item[1)] if $T \in \mathcal G_A$, then $0 < |T| < d |A|$;
\item[2)] every double coset of the group $G(i-1)$ over the pair of subgroups $\langle A \rangle$, $\langle A \rangle$ contains at most
one word from $\mathcal G_A$, and this word has minimal length among words from the coset.
\end{itemize}

The set of defining words $\mathcal S_i$ of rank $i$ is
constructed as follows. First, include in $\mathcal S_i$ all words
$A^{n_0}$ for each $A \in \mathcal X_i$. Furthermore, for every $A
\in \mathcal X_i$, if $a_j \notin \langle A \rangle \subset
G(i-1)$, and $a_j \notin \langle A \rangle a_k \langle A \rangle$
for $k<j$, then for every $T \in \mathcal G_A \setminus (\langle A
\rangle a_j \langle A \rangle)$ we include in $\mathcal S_i$ the
word
$$
a_jA^{n+j-1}TA^{n+m+j-1}\dots TA^{n+m(h-1)+j-1} ,
$$
where $j$ runs from $1$ to $m$.

Then set $\mathcal R_i =\mathcal R_{i-1} \cup \mathcal S_i$,
define $G(i)=\langle a_1,~a_2, \dots , a_m | R=1, ~R \in \mathcal
R_i \rangle$ and finally,
$$
G(\infty)=\langle a_1,~a_2, \dots ,a_m | R=1, ~R \in \mathcal R =
\bigcup_{i=0}^{\infty} \mathcal R_i \rangle
$$

The following lemmas can be found in \cite{Sem-diss} (see also
\cite[Chapter 8, 9]{Ols-book})


\begin{lem} \label{2gen}
The group $G(\infty)$ is infinite. Every proper subgroup of
$G(\infty)$ is cyclic of order dividing $n_0$.
\end{lem}

\begin{lem} \label{conj}
Any nontrivial element of $G(\infty)$ is conjugate to a power of a
period of certain rank. The centralizer of any nontrivial element
of $G(\infty)$ is cyclic of order $n_0$.
\end{lem}

For the proof of the following lemma we refer to \cite{Osin-unif},
Corollary 2.3 :
\begin{lem} \label{comm}
If $a,b \in G(\infty)$, $[a,~b] \ne 1$, then for every $i\in
\mathbb Z$, we have $[a^i,~b] \ne 1$ whenever $a^i \ne 1$.
\end{lem}

In the next two lemmas additional parameter $\sigma$ appears. It
corresponds to $100k_1^{-3/2}\zeta^{-1}$ from \cite{Sem-diss}.
\begin{lem} \label{est1}
Let $C$ be a period of a certain rank, $V \equiv C^k$, where
$\sigma < k <n_0/2$. Suppose that an element $W$ does not commute
with $V$ and has minimal length among the elements of $\langle C^k
\rangle W \langle C^k \rangle$. Assume also that $[C^k,~W]$ is
conjugated to $A^l$, where $A$ is a period of a certain rank and
$|l| \le n_0/2$. Then $|l| \le \sigma$ and the pair
$([C^k,~W],~C^k)$ is conjugate to the pair $(A^l,~B)$, where
$|B|<d|A|$.
\end{lem}

\begin{lem}\label{est2}
Let $l$ be an integer and $A,~B$ be elements of $G(\infty)$ such
that
\begin{enumerate}
\item $A$ is a period of certain rank;

\item $[A,~B] \ne 1$ in $G(\infty)$;

\item $|l| \le \sigma$ and $|B|<d|A|$.
\end{enumerate}
Then there exists an integer $s$ such that the pair $(BA^{ls},~B)$
is conjugate to a pair $(F,~T)$, where $F$ is a period of certain
rank, $[F,~T] \ne 1$ and $|T|<d|F|$.
\end{lem}

\section{Proof of the main results}


\begin{proof}[Proof of Theorem 1.1] We show that the group $G(\infty)$
construction of which was outlined in Section 2 has property (U)
with respect to the generating set $\{a_1,~a_2, \dots a_m\}$. In
view of Lemma \ref{2gen} it suffices to consider generating sets
$Y$ of $G(\infty)$ consisting of two non-commuting elements. First
we will consider the situation when $Y$ consists of a period of
some rank and a "short" word. Then it will be shown how to reduce
general situation to the considered one.

Assuming that generating set $Y$ consists of a period $F$ and a
word $T$ such that $[F,T] \ne 1$ and $|T| <d|F|$, we show that
there is a uniform bound (depending on the group $G(\infty)$ only)
for the lengths of elements $a_1, \dots , a_m$ as words in $F$ and
$T$. Denote by $i$ the rank of period $F$ and consider the
relations imposed on the $i$-th step.

If $a_j \in \langle F \rangle \subset G(i-1)$ or $T \in \langle F
\rangle a_j \langle F \rangle \subset G(i-1)$, then, in the group
$G(\infty)$, $|a_j|_{\{F,T\}}<n_0$ in view of the relation
$F^{n_0}=1$. Suppose now that $a_j \notin \langle F \rangle \cup
(\langle F \rangle T \langle F \rangle) \subset G(i-1)$. If $l$ is
minimal index such that $a_l \in \langle F \rangle a_j \langle F
\rangle$, then for some element $T_1 \in \langle F \rangle T
\langle F \rangle$ we encounter a relation
$$
a_lF^{n+j-1}T_1F^{n+m+j-1}\dots T_1F^{n+m(h-1)+j-1}=1~.
$$

Again, in view of $F^{n_0}=1$, one has $|T_1|_{\{F,T\}}<n_0$ and
therefore $|a_l|_{\{F,T\}}<n_0(n_0+m)$ in the group $G(\infty)$.
Consequently, $|a_j|_{\{F,T\}}<n_0(n_0+m+1)$ for any $a_j \in
\langle F \rangle a_l \langle F \rangle$ in $G(\infty)$. Thus,
every $a_j$ is equal to a word in $F$ and $T$ of length at most
$n_0(n_0+m+1)$ in the group $G(\infty)$.


Let now $Y$ consist of two non-commuting elements $u$ and $v$. We
will show that there exist a constant $M$ depending on the group
$G(\infty)$ only such that the following condition holds. There
exists a period $F$ and a word $T$, $[F,~T] \ne 1$, $|T|<d|F|$,
such that
$$
\max(|t^{-1}Ft|_{\{u,v\}}, |t^{-1}Tt|_{\{u,v\}})<M
$$
for some element $t \in G(\infty)$.

By Lemma \ref{conj}, we have $u=P^{-1}C^{k^{\prime}}P$ for some
element $P$ and a period $C$ of certain rank. Without loss of
generality we can assume that $0 < k^{\prime} < n_0/2$. Recall
that $\sigma <n_0/4$. There is a number $i \in \mathbb N$ such
that
$$
i < \frac{n_0}{2}
$$
and $\sigma < ik^{\prime} < n_0/2$. Thus we have
$$
u^i=P^{-1}C^kP~,
$$
where
$$
\sigma < k < \frac{n_0}{2}~.
$$
Using Lemma \ref{comm}, we note that
$$
[C^k,~PvP^{-1}]=P[u^i,~v]P^{-1} \ne 1.
$$
Denote by $W_0$ the element $PvP^{-1}$ and by $W$ the shortest
element in the double coset $\langle C^k \rangle W_0 \langle C^k
\rangle$. There are some integers $k_1,~k_2$ satisfying the
inequality $\max(|k_1|,~|k_2|) < n_0/2$, such that $W=C^{kk_1} W_0
C^{kk_2}$. Since $W=Pu^{ik_1}vu^{ik_2}P^{-1}$, we obtain the
following inequality:
$$
|P^{-1}WP|_{\{u,~v\}} \le 2i+k_1+k_2 +1 <2n_0.
$$
The elements $C^k$ and $W$ do not commute with each other.
Therefore, by Lemma \ref{conj}, $[C^k,~W]$ is conjugate to $A^l$,
where $A$ is a period of a certain rank:
$$
[C^k,~W]=Q^{-1}A^lQ
$$
for some element $Q$. Set
$$
B=QC^kQ^{-1}.
$$
Applying Lemma \ref{est1}, we obtain $|l| \le \sigma$ and $|B| <
d|A|$. Note that all conditions of Lemma \ref{est2} are satisfied
for $A$ and $B$. Therefore, there is an integer $s$ such that the
pair $(BA^{ls},~B)$ is conjugate to a pair $(F,~T)$, where $F$ is
a period of a certain rank, $[F,~T] \ne 1$, and $|T| <d|F|$.
We can assume that $s<n_0$, so that
$$
\max(|Q^{-1}BA^{ls}Q|_{\{[C^k,W],C^k\}},
|Q^{-1}BQ|_{\{[C^k,W],C^k\}})<n_0+1.
$$
Note that
$$
\max(|P^{-1}[C^k,W]P|_{\{u,v\}}, |P^{-1}C^kP|_{\{u,v\}})<5n_0.
$$
Therefore for some element $t$,
$$
\max(|t^{-1}Ft|_{\{u,v\}}, |t^{-1}Tt|_{\{u,v\}})<5n_0(n_0+1).
$$
To complete the proof it suffices to set $M=5n_0(n_0+1)$ and
$C=C(G(\infty))=Mn_0(n_0+m+1)$.
\end{proof}

To prove Corollary 1.2 we need two elementary facts about Kazhdan
constants.

\begin{lem}\label{GH}
Let $G$ be a group, $X$ a finite generating set of $G$. Then:
\begin{enumerate}
\item $\varkappa (G,X)=\varkappa (G, t^{-1}Xt) $ for any element
$t\in G$.

\item For any finite generating set $Y$ in $G$, $\varkappa (G,Y)
\ge \varkappa (G,X)/d$, where $d=\max\limits_{x\in X} |x|_Y$.
\end{enumerate}
\end{lem}

\begin{proof} For any $g,t\in G$, any unitary presentation $\pi
:G\to U(\mathcal H)$, and any unit vector $\xi \in \mathcal H$, we
have
$$ \| \pi (t^{-1}gt)\xi -\xi \| =\| \pi (gt)\xi -\pi (t)\xi \| =
\| \pi (g)\xi -\xi \| .$$ This implies the first assertion of
the lemma.

Further let $x=y_1y_2\cdots y_l $ for some $y_1, \ldots , y_l\in
Y\cup Y^{-1} $ and let $\xi \in \mathcal H^1$. Then
$$
\begin{array}{rl}
\| \pi(x)\xi - \xi\|  \le & \|\pi(y_1)\xi - \xi \| +
\sum\limits_{i=1}^{l-1} \| \pi(y_1y_2 \cdots y_{i+1})
\xi-\pi(y_1y_2\cdots y_{i}) \xi\| = \\&\\
& \sum\limits_{i=1}^{l} \|\pi(y_i)\xi - \xi \| \le l
\max\limits_{y \in Y} \| \pi(y)\xi - \xi \|.
\end{array}
$$
This yields the second assertion.
\end{proof}

\begin{proof}[Proof of Corollary \ref{ut}]
Applying Theorem \ref{U} to a torsion free hyperbolic Kazhdan
group, one obtains group $G(\infty)$ satisfying both properties
(T) and (U). We denote by $\varkappa=\varkappa(G(\infty), X)$ the
Kazhdan constant of $G(\infty)$ with respect to the generating set
$X=\{a_1,~a_2, \dots a_m\}$. If $Y$ is any other generating set of
$G(\infty)$, then, by Lemma \ref{GH},
$$
\varkappa(G(\infty), Y) \ge \frac1C \varkappa(G(\infty), X)~,
$$
where $C=C(G(\infty))$ is the constant realizing property (U) of
$G(\infty)$ with respect to $X$. Consequently, property (T) of
$G(\infty)$ implies $\varkappa(G(\infty))>0$.
\end{proof}



\begin{thebibliography}{99}

\bibitem{BHV}
B. Bekka, P. de la Harpe, A. Valette, \emph{Kazhdan's Property
(T)}; available at
http://www.unige.ch/math/biblio/preprint/2006/KazhdansPropertyT.pdf

\bibitem{GZ}
T.Gelander, A. \. Zuk, \emph{Dependence of Kazhdan constants on
generating subsets}, Israel Journal of Math., vol. 129, 2002, pp.
93--99.

\bibitem{Gro87} M. Gromov, \emph{Hyperbolic groups}, Essays in
Group Theory (S.M.Gersten, editor), Math. Sci. Res. Inst. Publ.,
vol 8, Springer-Verlag, Berlin, 1987, pp.  75--263.





\bibitem{IO}
S.V. Ivanov, A.~Yu. Ol'shanskii, \emph{Hyperbolic groups and their
quotients of finite exponent}, Trans. Amer. Math. Soc. 348 (1996)
N 6, 2091--2138.


\bibitem{Kazhd}
D.A. Kazhdan, \emph{On the connection of the dual space of a group
with the structure of its closed subgroups (Russian)}, Funkcional.
Anal. i Prilo\v zen., {\bf 1} (1967), pp. 71--74.

\bibitem{Lub}
A. Lubotzky, \emph{Discrete groups, expanding graphs and invariant
measures}, vol. 125 of Progress in Mathematics, Birkhauser Verlag,
Basel;Boston;Berlin, 1994.

\bibitem{Ols_cycl}
A.~Yu. Ol'shanskii, \emph{Infinite groups with cyclic subgroups},
Doklady Akad. Nauk SSSR, 245 (1979), N 4, pp. 785-787 (in
Russian); English translation: Soviet Math. Dokl. 20 (1979), no.
2, pp. 343--346.

\bibitem{Ols}
A.~Yu. Ol'shanskii, \emph{Periodic factor groups of hyperbolic
groups}, (Russian) Mat. Sbornik {\bf 182} (1991), no. 4;
translation in Math. USSR Sbornik Vol. 72 (1992), no. 2, pp.
519--541.

\bibitem{Ols-book}
A.~Yu. Ol'shanskii, \emph{The geometry of defining relations in groups}, Nauka, Moscow, 1989; Translated in Math and Its Applications
(Soviet series), \textbf{70}, Kluwer Acad. Publishers, 1991.

\bibitem{Osin-unif}
D.V. Osin, \emph{ Uniform non--amenability of free Burnside
groups}, submitted to Mat. Zametki (Math. Notes); available at
http://www.arxiv.org/abs/math.GR/0404073

\bibitem{Osin-kazh}
D.V. Osin, \emph{Kazhdan constants of hyperbolic groups,}
(Russian) Funktsional. Anal. i Prilozhen. {\bf 36} (2002), no. 4,
pp. 46--54; translation in Funct. Anal. Appl. {\bf 36} (2002), no.
4, pp. 290--297.

\bibitem{Sem-diss}
Yu.S. Semenov, \emph{Some quotient groups and rings of hyperbolic
groups}, PhD Thesis, Moscow State University, 1994.

\bibitem{Sem-obz}
Yu.S. Semenov, \emph{Some quotient groups of hyperbolic groups,}
(Russian) Vestnik Moskov. Univ. Ser. I Mat. Mekh. 1993, no. 3, pp.
88--90; translation in Moscow Univ. Math. Bull. 48 (1993), no. 3,
pp. 39--40.

\bibitem{Zelm1}
E.I.Zelmanov, \emph{Solution of the restricted Burnside problem
for groups of odd exponent,} Math. USSR Izv {\bf 36} (1991), pp.
41--60.


\bibitem{Zuk}
A. \. Zuk, \emph{Property (T) and Kazhdan constants for discrete
groups,} GAFA, Geom. Funct. Anal. Vol. 13 (2003) pp. 643--670.

\end{thebibliography}
\end{document}